\input epsf

\magnification=\magstep1

\font\titre=cmbx10 scaled\magstep2
\font\soustitre=cmbx10 scaled\magstep1
\font\aut=cmr10 scaled\magstep1

\def\square{{\vcenter{\hrule height.4pt
      \hbox{\vrule width.4pt height5pt \hskip5pt
	   \vrule width.4pt}
      \hrule height.4pt}}}
\def\qed{\hfill$\square$}
\def\semidirect{\mathop
 {\rlap{\kern 7pt\vrule height 5pt width .4pt depth 0pt} \times}}

\def\HH{{\cal H}}
\def\MM{{\cal M}}
\def\BZ{{\bf Z}}
\def\SL{{\rm SL}}
\def\Stab{{\rm Stab}}
\def\inf{{\rm inf}}
\def\TT{{\cal T}}
\def\OO{{\cal O}}
\def\Ind{{\rm Ind}}
\def\fd{{\rm fd}}
\def\PMM{{\cal PM}}
\def\Cub{{\rm Cub}}
\def\BR{{\bf R}}
\def\Ker{{\rm Ker}}
\def\Im{{\rm Im}}
\def\Com{{\rm Com}}

\phantom{blabla}
\bigskip\bigskip
\centerline{\titre Actions and irreducible representations}
\medskip
\centerline{\titre of the mapping
class group}

\bigskip
\centerline{\aut Luis Paris}

\bigskip
\centerline{\aut May 1999}

\bigskip\bigskip
\centerline{\bf Abstract}

\bigskip
Let $G$ be a countable discrete group. Call two subgroups
$H_1$ and $H_2$ of $G$ commensurable if $H_1\cap H_2$ has
finite index in both $H_1$ and $H_2$. We say that an
action of $G$ on a discrete set $X$ has noncommensurable
stabilizers if the stabilizers of any two distinct points
of $X$ are not commensurable. We prove in this paper that
the action of the mapping class group on the complex of
curves has noncommensurable stabilizers. Following a
method due to Burger and de la Harpe, this action
leads to constructions of irreducible unitary
representations of the mapping class group.

\bigskip\noindent
{\it Mathematics Subject Classifications:} Primary 57N05;
Secondary 20F38, 22D10, 22D30.

\bigskip\bigskip
\centerline{\soustitre 1. Introduction}

\bigskip
Throughout the paper, $M$ will denote a compact connected oriented
surface, possibly with boundary, and $P=\{p_1,\dots,p_m\}$
a finite collection of $m$ points (called {\it punctures})
in the interior of $M$. Define $\HH(M,P)$ to be the group
of orientation-preserving homeomorphisms $h:M\to M$ such
that $h$ is the identity on the boundary of $M$ and
$h(P)=P$. Two homeomorphisms $h$ and $h'$ in $\HH(M,P)$
are said to be {\it isotopic} if there exists a continuous
family $h_t\in\HH(M,P)$, $t\in [0,1]$ such that $h_0=h$
and $h_1=h'$. The {\it mapping class group} $\MM(M,P)$ is
the group of isotopy classes in
$\HH(M,P)$. 

Let $G$ be a countable discrete group. An {\it unitary}
{\it representation} is a homomorphism $\pi$ from $G$ to
the unitary group of a separable Hilbert space $V$. We say
that $\pi$ is {\it irreducible} if $V$ does not contain
any closed subspace invariant by the action of $G$,
except $\{0\}$ and $V$ itself. The space of equivalence
classes of irreducible unitary representations of $G$ is called
the {\it unitary} {\it dual} of $G$ and is denoted by
$\hat G$. We will also denote by $\hat G^\fd$ the subspace
of $\hat G$ of equivalence classes of finite dimensional
irreducible unitary representations. Besides the particular case
of virtually abelian groups (see [War, Theorem 5.4.1.4]),
there is no systematic
procedure of constructing all irreducible unitary representations
of $G$. However, one may hope to construct large families
of them.

Recall that two subgroups $H_1$ and $H_2$ of $G$
are {\it commensurable} if $H_1\cap H_2$ has finite
index in both $H_1$ and $H_2$. We say that
an action of $G$ on a discrete set $X$ has 
{\it noncommensurable stabilizers} if the stabilizers of
any two distinct points of $X$ are not commensurable.
A method to construct
irreducible unitary representations
is proposed by Burger and
de la Harpe in [BH]. This method consists first in
translating some results due to Mackey [Mac] into the
terminology of the discrete group theory, and then
to apply these Mackey's results to the situation of an
action with noncommensurable stabilizers.

With a punctured surface $(M,P)$ one can associate a
simplicial complex $X(M,P)$ called {\it complex of
curves}. This complex has been introduced by Harvey [Harv]
and plays a prominent r\^ole in the study of mapping class
groups (see [BLM], [Ha1], [Ha2], [Iv1], [Iv2], [IM],
[McC], [MM1], [MM2]). 

We consider in this paper $X(M,P)$ as a set, ignoring the
simplicial structure. Our main result is that the action
of $\MM(M,P)$ on $X(M,P)$ has noncommensurable
stabilizers. We also describe the structure of the
stabilizers, leading by this way to constructions of
irreducible unitary representations of $\MM(M,P)$. Furthermore,
this action provides interesting phenomena that show the
limits of a possible generalization of Mackey's results.

Our work is organized as follows. In Section 2 we give the
definition of the complex of curves $X(M,P)$ and study 
some basic examples. In Setion 3 we state the method of
Burger and de la Harpe of constructing irreducible
unitary representations. In Section 4 we state some preliminary
results on curves and Dehn twists which will be useful
for the remainder of the paper. These preliminary results
are enterely proved in [PR2]. Section 5 is dedicated
to our main result and its consequences, and we study in
Section 6 the structure of the stabilizers.

\bigskip\bigskip
\centerline{\soustitre 2. Definitions and examples}

\bigskip
By a {\it simple closed curve} in $M\setminus P$ we mean
an embedding $a:S^1\to M\setminus P$. Note that $a$ has an
orientation; the curve with opposite orientation but same
image will be denoted by $a^{-1}$. By abuse of notation,
we also use $a$ for the image of $a$. We say that $a$ is
{\it generic} if it does not bound a disk in $M$
containing zero or one puncture. Two simple closed curves
$a,b:S^1\to M\setminus P$ are {\it isotopic} if there
exists a continuous family $a_t:S^1\to M\setminus P$,
$t\in [0,1]$ of simple closed curves such that $a_0=a$ and
$a_1=b$. Isotopy of curves is an equivalence relation
which we denote by $a\simeq b$. Note that $h(a)\simeq
h'(b)$ if $a\simeq b$ and if $h,h'\in\HH(M,P)$ are isotopic.
So, the mapping class group $\MM(M,P)$ acts on the set of
isotopy classes of simple closed curves.

Define a {\it generic $r$-family of disjoint curves} to be a 
$r$-tuple $(a_1,\dots,a_r)$ of generic closed curves such
that:

\smallskip
(a) $a_i\cap a_j=\emptyset$, for $i\neq j$;

\smallskip
(b) $a_i$ is neither isotopic to $a_j$ nor $a_j^{-1}$, for
$i\neq j$;

\smallskip
(c) $a_i$ is not isotopic with a boundary curve of $M$.

\smallskip\noindent
We say that two generic $r$-families of disjoint curves
$A=(a_1,\dots,a_r)$ and $B=(b_1,\dots,b_r)$ are {\it
equivalent} if there exists a permutation
$\sigma\in\Sigma_r$ such that $a_i\simeq
b_{\sigma(i)}^{\pm 1}$ for each $i=1,\dots, r$. 
We will write $[a_1,\dots,a_r]$ for the equivalence class
of $(a_1,\dots,a_r)$.
Then
$X_r(M,P)$ denotes the set of equivalence classes of
generic $r$-families of disjoint curves. 
Notice that $X_1(M,P)$ is the set of isotopy classes
of nonoriented generic curves.
The action of $\MM(M,P)$
on $X_r(M,P)$ is induced by the action of $\MM(M,P)$ on
the set of isotopy classes of simple closed curves. Now, define
the {\it complex of curves} to be:
$$
X(M,P)=\bigcup_{r=1}^\infty X_r(M,P)\ .
$$

\noindent
{\bf Example 1.} $\MM(S^2)=\MM(S^2,\{p\})=
\MM(D^2)=\MM(D^2,\{p\})=\{1\}$. 
These are
the only mapping class groups which are trivial.
$\MM(S^2,\{p_1,p_2\})$ is a cyclic group of order two, and
$\MM(S^2,\{p_1,p_2,p_3\})$ is isomorphic with the group 
$\Sigma_3$ of
permutations of $\{1,2,3\}$. These are the only 
non-trivial mapping class groups which are finite. Note
that, in these examples, $X(M,P)$ is empty. However, the
group being finite, the irreducible unitary representations of
$\MM(M,P)$ are well understood.

\bigskip\noindent
{\bf Example 2.} There are four abelian mapping class
groups. The first one is the mapping class group of the
cylinder: $\MM(S^1\times I)=\BZ$. The other three cases
are of importance for the remainder of the paper. 
Choose two points $p_1,p_2$ in the interior of the
disk $D^2$ and call a pair equivalent to $(D^2,\{p_1,p_2\})$
a {\it pantalon of type I}. Choose a point $p$ in the interior
of the cylinder $S^1\times I$ and call a pair equivalent
to $(S^1\times I,\{p\})$ a {\it pantalon of type II}.
Call a pair equivalent to $(M,\emptyset)$, where
$M$ is the sphere with three holes, a {\it pantalon of
type III}.  
$\MM(D^2,\{p_1,p_2\})$ is the infinite cyclic group
generated by the ``half-twist'' which interchanges the
punctures. $\MM(S^1\times I,\{p\})$ is isomorphic with $\BZ^2$
and is generated by the Dehn twists along the boundary
components. Let $M$ be the sphere with three holes.
$\MM(M)$ is isomorphic with $\BZ^3$ and is generated by
the Dehn twists along the boundary components. Like in the
finite case, $X(M,P)$ is empty in these four examples.
Nevertheless, the group
being abelian, the irreducible representations of
$\MM(M,P)$ are well understood. Note that there is no
other example where $X(M,P)$ is empty. Moreover, if
$\MM(M,P)$ is virtually abelian, then $(M,P)$ is as
in Example 1 or 2, thus either $\MM(M,P)$ is
finite or free abelian.   

\bigskip\noindent
{\bf Example 3.} The mapping class group of the torus
$T^2=S^1\times S^1$ with either zero or one punture is the
modular group of $2\times 2$ matrices with integer
coefficients: $\MM(T^2)=\MM(T^2,\{p\})=\SL(2,\BZ)$. There is
only one $\MM(T^2)$-orbit in $X(T^2)$ represented by
$[a]$, where $a=S^1\times\ast$. The stabilizer of $[a]$ is the
subgroup of upper triangular matrices in $\SL(2,\BZ)$
(with diagonal entries $\pm 1$) which is known to be its
own commensurator in $\SL(2,\BZ)$ by [BH].

\bigskip
\epsfysize=2.5 truecm
\centerline{\epsfbox{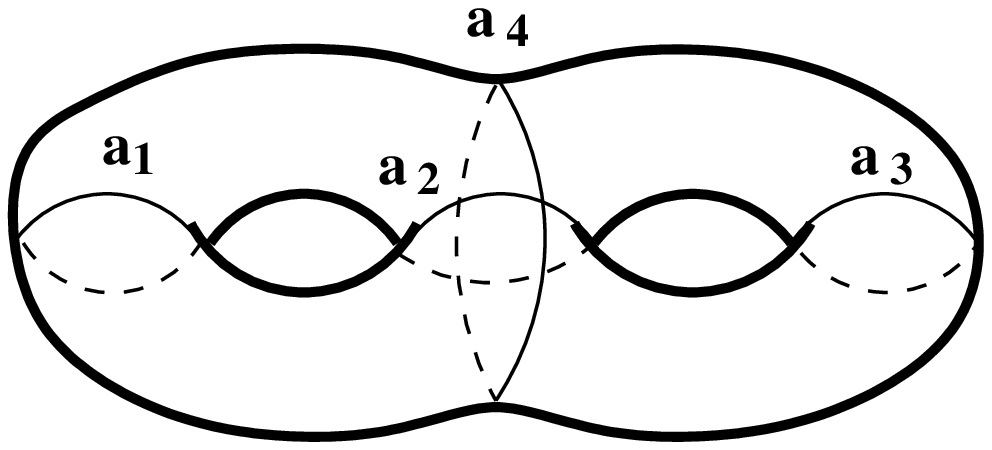}}
\bigskip
\centerline{Figure 1: Some curves in the closed surface of 
genus $2$}

\bigskip\noindent
{\bf Example 4.} Assume that $M$ is a closed surface of genus
$2$, and let $a_1,a_2,a_3,a_4$ be the simple closed curves
drawn in Figure 1. There are $6$ $\MM(M)$-orbits in $X(M)$
represented by $[a_1]$, $[a_4]$, $[a_1,a_2]$, $[a_1,a_4]$,
$[a_1,a_2,a_3]$, $[a_1,a_3,a_4]$.

\bigskip\noindent
{\bf Example 5.} The mapping class group of the disk with
$m=|P|$ punctures is known to be the Artin braid group
$B_m$ on $m$ strings. Consider a simple closed curve
$a:S^1\to D^2\setminus P$. Then $a$ bounds a disk $D_a$
embedded in $D^2$. The group $B_a=\MM(D_a,P\cap D_a)$ is a
copy of the Artin braid group  on $|P\cap D_a|$ strings
which can be assumed to be embedded in $B_m=\MM(D^2,P)$
via the homomorphism $\MM(D_a,P\cap D_a)\to\MM(D^2,P)$
induced by the inclusion $D_a\subset D^2$. 
It is a non-trivial but well known fact that the above 
homomorphism is injective.
Such a subgroup
is called {\it geometric subgroup} of $B_m$. Assume $a$ is
generic. Rolfsen proved in [Ro] that the normalizer of
$B_a$ in $B_m$ is equal to its commensurator in $B_m$, and
is equal to the stabilizer of $[a]$. He also proved that the action
of $B_m$ in $X_1(D^2,P)$ has noncommensurable stabilizers, 
infering by this way
that the stabilizer of $[a]$ is its own commensurator in $B_m$. Some
of these results have been generalized to parabolic
subgroups of Artin groups [Pa], to geometric subgroups of
surface braid groups [PR1], and to geometric subgroups of
mapping class groups [PR2].

\bigskip\bigskip
\centerline{\soustitre 3. Constructing irreducible 
representations}

\bigskip
Let $G$ be a countable discrete group.
Define the {\it commensurator} of a
subgroup $H$ of $G$ to be:
$$
\Com_G(H)=\{g\in G\ ;\ gHg^{-1}\ {\rm and}\ H\ {\rm are\
commensurable}\}\ .
$$

\noindent
{\bf Theorem 3.1} (Mackey [Mac], Burger, de la Harpe [BH]).
{\sl Let $H$ be a subgroup of $G$. The following
statements are equivalent:

\smallskip
(1) $\Com_G(H)=H$;

\smallskip
(2) the left regular representation $\lambda_{G/H}$ of $G$
in $l^2(G/H)$ is irreducible;

\smallskip
(3) for all finite dimensional irreducible unitary representation
$\pi$ of $H$, the induced representation $\Ind_H^G\pi$ is
irreducible.

\smallskip\noindent
Furthermore, if (1)-(3) are satisfied (i.e. $\Com_G(H)=H$),
then unitary induction provides a well defined and
injective map
$$
\Ind_H^G:\hat H^\fd\to\hat G\ .\quad\square
$$}

\noindent
{\bf Theorem 3.2} (Mackey [Mac], Burger, de la Harpe [BH]).
{\sl Let $H_1$ and $H_2$ be two subgroups of $G$ such that
$\Com_G(H_1)=H_1$ and $\Com_G(H_2)=H_2$. The following
statements are equivalent:

\smallskip
(1) $H_1$ and $H_2$ are conjugate;

\smallskip
(2) the left regular representations $\lambda_{G/H_1}$ and
$\lambda_{G/H_2}$ are equivalent;

\smallskip
(3) there exist finite dimensional irreducible unitary
representations $\pi_1$ and $\pi_2$ of $H_1$ and $H_2$,
repectively, such that $\Ind_{H_1}^G\pi_1$ and
$\Ind_{H_2}^G\pi_2$ are equivalent.

\smallskip\noindent
In particular, in case $H_1$ and $H_2$ are not
conjugate, then unitary induction provides a well
defined and injective map
$$
\Ind:\hat H_1^\fd\amalg\hat H_2^\fd\to\hat G\ .\quad\square
$$}

Let $G\times X\to X$ be an action of $G$ on some discrete
set $X$. For $x\in X$, let
$$
\Stab(x)=\{g\in G\ ;\ gx=x\}
$$
denote the {\it stabilizer} of $x$. We say that the action
has {\it noncommensurable stabilizers} if $\Stab(x)$ and
$\Stab(y)$ are not commensurable for all $x,y\in X$,
$x\neq y$.

\bigskip\noindent
{\bf Proposition 3.3} (Burger, de la Harpe [BH]). {\sl Let
$G\times X\to X$ be an action with
noncommensurable stabilizers.

\smallskip
(i) $\Com_G(\Stab(x))=\Stab(x)$ for all $x\in X$;

\smallskip
(ii) $\Stab(x)$ and $\Stab(y)$ are conjugate if and
only if $x$ and $y$ are in the same orbit.\qed}

\bigskip\noindent
{\bf Corollary 3.4} (Burger, de la Harpe [BH]). {\sl Let
$G\times X\to X$ be an action with
noncommensurable stabilizers.
Let $\OO$ denote the set of orbits in $X$. Then
unitary induction provides a well defined and
injective map
$$
\Ind: \coprod_{x\in \OO}\widehat{\Stab(x)}^\fd\to\hat G\ .\quad\square
$$}

\bigskip\bigskip
\centerline{\soustitre 4. Curves and Dehn twists}

\bigskip
As pointed out in the introduction, the statements
of this section are proved in [PR2]. So, we do
not give the proofs here and refer to [PR2] for them.

The {\it index of intersection} of two simple closed curves
$a,b:S^1\to M\setminus P$ is:
$$
I(a,b)=\inf\{|a'\cap b'|\ ;\ a'\simeq a\ {\rm and}\ b'
\simeq b\}\ .
$$
Note that:

\smallskip
(1) if $a$ is not generic, then $I(a,b)=0$ for every simple
closed curve $b$;

\smallskip
(2) if $a\simeq b$, then $I(a,b)=0$.

\bigskip\noindent
{\bf Proposition 4.1.} {\sl Let $(a_1,\dots, a_r)$ be a
generic family of disjoint curves. For every $i\in\{1,\dots, r\}$,
there exists a simple closed curve $b:S^1\to M\setminus P$
such that $a_j\cap b=\emptyset$ for $j\neq i$, and
$|a_i\cap b|=I(a_i,b)>0$.}\qed

\bigskip\noindent
{\bf Proposition 4.2.} {\sl Let $(a_1,\dots,a_r)$ be a
generic family of disjoint curves and $b:S^1\to M\setminus P$ a
simple closed curve such that $I(a_i,b)=0$ for all
$i\in\{1,\dots,r\}$. Then there exists a simple closed
curve $b':S^1\to M\setminus P$ isotopic with $b$ and such
that $a_i\cap b'=\emptyset$ for all $i\in\{1,\dots,r\}$.}
\qed

\bigskip\noindent
{\bf Proposition 4.3.} {\sl Let $(a_1,\dots, a_r)$ and
$(b_1,\dots, b_r)$ be two generic $r$-families of disjoint curves
such that $a_i\simeq b_i$ for all $i\in\{1,\dots,r\}$.
Then there exists an isotopy $h_t\in\HH(M,P)$,
$t\in[0,1]$, such that $h_0=$identity and $h_1\circ
a_i=b_i$ for all $i\in\{1,\dots, r\}$.} \qed

\bigskip
By a {\it subsurface} $N$ of $M$ we mean a closed subset
which is also a surface and for which we always assume
$\partial N\cap P=\emptyset$. We say furthermore that $N$
is {\it essential} if each component of
$\overline{M\setminus N}$ which is a disk has a nonempty
intersection with the puncture set $P$.

\bigskip\noindent
{\bf Proposition 4.4.} {\sl Let $N$ be an essential
subsurface of $M$, and let $a,b:S^1\to N\setminus N\cap P$
be two generic closed curves. Assume that $a$ is not
isotopic with a boundary curve of $N$. Then $a$ is isotopic
with $b$ in $M\setminus P$ if and only if $a$ is isotopic
with $b$ in $N\setminus N\cap P$.}\qed

\bigskip
Let $a:S^1\to M\setminus P$ be a simple closed curve
interior to $M$.  Choose a tubular neighborhood $v$ of $a$
(namely, an oriented embedding $v:S^1\times I\to
M\setminus P$ such that $v(z,1/2)=a(z)$ for all $z\in
S^1$), and define the homeomorphism $T_a\in\HH(M,P)$ by:
$$
(T_a\circ v)(z,t)=v(e^{2i\pi t}z,t)\ ,
$$
and $T_a$ is the identity on the exterior of the image of
$v$ in $M$. The {\it Dehn twist} along $a$ is defined to
be the element $\tau_a$ of $\MM(M,P)$ represented by
$T_a$. Note that:

\smallskip
(1) the definition of $\tau_a$ does not depend on the choice
of the tubular neighborhood;

\smallskip
(2) $\tau_a=\tau_b$ if $a\simeq b$;

\smallskip
(3) $\tau_a=\tau_{a^{-1}}$;

\smallskip
(4) $\tau_a=1$ if $a$ is not generic;

\smallskip
(5) $T_a(a)=a$;

\smallskip
(6) let $\xi\in\MM(M,P)$, and let $h\in\HH(M,P)$ which represents
$\xi$; then $\xi\tau_a\xi^{-1}$ is the Dehn twist along
$h(a)$;

\smallskip
(7) $\tau_a$ and $\tau_b$ commute if $I(a,b)=0$.

\smallskip
The following proposition is a special case of a formula in
[FLP].

\bigskip\noindent
{\bf Proposition 4.5.} {\sl Let $a,b:S^1\to M\setminus P$
be two simple closed curves and $n\in\BZ$. Then
$$
I(\tau_a^n(b),b)=|n|I(a,b)^2\ .\quad\square
$$}

The next two propositions can be also found in [IM].

\bigskip\noindent
{\bf Proposition 4.6.} {\sl Let $a,b:S^1\to M\setminus P$
be two generic closed curves and
$j,k\in\BZ\setminus\{0\}$. If $\tau_a^j=\tau_b^k$, then
$a\simeq b$.} \qed

\bigskip\noindent
{\bf Proposition 4.7.} {\sl Let $a,b:S^1\to M\setminus P$ be two
generic closed curves and $j,k\in\BZ\setminus\{0\}$. If
$\tau_a^j$ and $\tau_b^k$ commute, then $I(a,b)=0$.}\qed

\bigskip\noindent
{\bf Proposition 4.8.} {\sl 
Consider $r$ generic closed curves $a_1,\dots, a_r$ such
that:

\smallskip
(a) $a_i\cap a_j=\emptyset$, for $i\neq j$;

\smallskip
(b) $a_i$ is neither isotopic to $a_j$ nor $a_j^{-1}$, for
$i\neq j$.

\smallskip\noindent
Then the subgroup
of $\MM(M,P)$ generated by
$\tau_{a_1},\dots,\tau_{a_r}$ is a free abelian group of
rank~$r$.} \qed 

\bigskip\noindent
{\bf Remark.} The $r$-tuple $(a_1,\dots, a_r)$ of the above
proposition in not necessarily a generic family of disjoint curves
in the sence that $a_i$ may be isotopic with a boundary
curve.

\bigskip\bigskip
\centerline{\soustitre 5. The action} 
 
\bigskip
Let $\alpha,\beta\in X(M,P)$. According to the terminology
on simplicial complexes, we say that $\alpha$ is a {\it
face} of $\beta$ if there exists a generic family of disjoint
curves $(a_1,\dots,a_r)$ such that
$\beta=[a_1,\dots,a_r]$ and $\alpha=[a_1,\dots, a_s]$
for some $s\in\{1,\dots,r\}$.

\bigskip\noindent
{\bf Proposition 5.1.} {\sl Let $\alpha,\beta\in X(M,P)$ such
that $\alpha$ is not a face of $\beta$. Then
the $\Stab(\beta)$-orbit of $\alpha$ is infinite.}

\bigskip\noindent
{\bf Proof.} Write $\alpha=[a_1,\dots,a_r]$ and
$\beta=[b_1,\dots,b_s]$. The hypothesis ``$\alpha$ is not
a face of $\beta$'' means that there is some $a_i$ which is
neither isotopic with $b_j$ nor $b_j^{-1}$, for all
$j\in\{1,\dots,s\}$.

Assume first that there exists some $b_j$ such that
$I(a_i,b_j)>0$. Consider the Dehn twist $\tau_{b_j}$ along
$b_j$.
Let $n\in\BZ\setminus \{0\}$. If
$\tau_{b_j}^n(\alpha)=\alpha$, then
$\tau_{b_j}^n(a_i)\simeq a_k^{\pm 1}$ for some $k\in\{1,\dots,
r\}$, thus
$$
I(\tau_{b_j}^n(a_i),a_i)=I(a_k^{\pm 1},a_i)=0\ .
$$
On the other hand, Proposition 4.5 implies
$$
I(\tau_{b_j}^n(a_i),a_i)=|n|I(a_i,b_j)^2\neq 0\ .
$$
So, $\tau_{b_j}^n(\alpha)\neq\alpha$ for all
$n\in\BZ\setminus\{0\}$. Since
$\tau_{b_j}\in\Stab(\beta)$, it follows that the
$\Stab(\beta)$-orbit of $\alpha$ is infinite.

Assume now that $I(a_i,b_j)=0$ for all $j\in\{1,\dots,s\}$.
By Proposition 4.2, we may assume that $a_i\cap
b_j=\emptyset$ for all $j\in\{1,\dots,s\}$, namely that
$(b_1,\dots,b_s,a_i)$ is a generic family of disjoint curves.
Proposition 4.1 provides a generic closed curve $c$ such
that $b_j\cap c=\emptyset$ for all $j\in\{1,\dots,s\}$ and
$|a_i\cap c|=I(a_i,c)>0$. As before,
$\tau_c^n(\alpha)=\alpha$ would imply
$I(\tau_c^n(a_i),a_i)=0$, but Proposition 4.5 guarantees
that $I(\tau_c^n(a_i),a_i)\neq 0$, thus
$\tau_c^n(\alpha)\neq\alpha$ for all
$n\in\BZ\setminus\{0\}$. Since $\tau_c\in\Stab(\beta)$, it
follows that the $\Stab(\beta)$-orbit of $\alpha$ is
infinite. \qed

\bigskip
An action of a group $G$ on a set $X$ is called {\it large}
if, for all $x,y\in X$, $x\neq y$, the $\Stab(x)$-orbit
of $y$ is infinite. It is shown in [BH] that large actions
have noncommensurable stabilizers.
Let $\alpha,\beta\in X_r(M,P)$, $\alpha\neq\beta$. 
Then $\alpha$
cannot be a face of $\beta$. So:

\bigskip\noindent
{\bf Corollary 5.2.} {\sl The action of $\MM(M,P)$ on
$X_r(M,P)$ is large for all $r$.} \qed

\bigskip
More generally:

\bigskip\noindent
{\bf Corollary 5.3.} {\sl The action of $\MM(M,P)$ on
$X(M,P)$ has noncommensurable stabilizers.}

\bigskip\noindent
{\bf Proof.} Let $\alpha,\beta\in
X(M,P)$, $\alpha\neq\beta$. We may assume that $\alpha$ 
is not a face of
$\beta$. Proposition 5.1 provides an infinite sequence
$\{g_i\}_{i\ge 0}$ in $\Stab(\beta)$ such that
$g_i\alpha\neq g_j\alpha$ for $i\neq j$. Then
$$
g_i(\Stab(\alpha)\cap\Stab(\beta))\neq g_j(\Stab(\alpha)
\cap\Stab(\beta))
$$
for $i\neq j$, hence $\Stab(\alpha)\cap\Stab(\beta)$ has
infinite index in $\Stab(\beta)$. \qed

\bigskip
So, by Proposition 3.3:

\bigskip\noindent
{\bf Corollary 5.4.} {\sl Let $\alpha\in X(M,P)$. Then
$\Stab(\alpha)$ is its own commensurator in $\MM(M,P)$.}
\qed

\bigskip
We turn now to the consequences of Proposition 5.1 on
irreducible unitary representations of the mapping class group. For
$\alpha\in X(M,P)$, we denote by $\lambda_\alpha$ the left
regular representation of $\MM(M,P)$ on
$l^2(\MM(M,P)/\Stab(\alpha))$.

\bigskip\noindent
{\bf Corollary 5.5.} {\sl (i) Let $\alpha\in X(M,P)$. Then
$\lambda_\alpha$ is an irreducible unitary representation
of $\MM(M,P)$.

\smallskip
(ii) Let $\alpha,\beta\in X(M,P)$. The representations
$\lambda_\alpha$ and $\lambda_\beta$ are equivalent if and
only if $\alpha$ and $\beta$ are in the same 
$\MM(M,P)$-orbit. In particular, $\lambda_\alpha$ and
$\lambda_\beta$ may not be equivalent if they have
different ranks.} \qed

\bigskip
Let $\OO_r$ denote the set of $\MM(M,P)$-orbits in
$X_r(M,P)$.

\bigskip\noindent
{\bf Corollary 5.6.} {\sl Unitary induction provides a well
defined and injective map
$$
\Ind:\coprod_{r=1}^\infty\left(\coprod_{\alpha\in\OO_r}
\widehat{\Stab(\alpha)}^{\fd}\right)\to\widehat{\MM(M,P)}
\ .\quad\square
$$}

\bigskip
A natural problem which arises from Corollaries 5.5 and 5.6
is to find a method to determine whether two classes
$\alpha,\beta\in X_r(M,P)$ are in the same 
$\MM(M,P)$-orbit, and to find a representative for each
$\MM(M,P)$-orbit. We see now that this question has a
natural answer in topological terms.

Let $A=(a_1,\dots,a_r)$ be a generic $r$-family of disjoint curves.
We denote by $M_A$ the natural compactification of
$M\setminus(\cup_{i=1}^ra_i)$ and by $\rho_A:M_A\to M$ the
continuous map induced by the inclusion map
$M\setminus(\cup_{i=1}^ra_i)\to M$. By abuse of notation,
we use $P$ for $\rho_A^{-1}(P)$. Let $N$ be a component of
$M_A$ and $c:S^1\to\partial N$ a boundary curve of $N$. We
say that $c$ is an {\it exterior boundary curve} of $N$ if
$\rho_A\circ c$ is a boundary curve of $M$. Let
$a_i:S^1\to M\setminus P$ be one of the curves of $A$. There
are two distinct boundary curves $c_i,c_i':S^1\to\partial
M_A$ of $M_A$ such that $\rho_A\circ c_i=\rho_A\circ
c_i'=a_i$. Let $N$ and $N'$ be the components of $M_A$
having $c_i$ and $c_i'$ as boundary curves, respectively.
We say that $a_i$ is a {\it separating limit curve} of $N$
if $N\neq N'$, and $a_i$ is a {\it nonseparating limit
curve} of $N$ if $N=N'$.

We say that $A$ determines a {\it pantalon decomposition}
(or simply $A$ is a {\it pantalon decomposition}) if, for each
component $N$ of $M_A$, $(N,P\cap N)$ is a pantalon (of
any type). Such a family exists if and only if neither
$(M,P)$ is as in Example 1 or 2, nor
$(M,P)=(T^2,\emptyset)$. A generic family
$A=(a_1,\dots,a_r)$ is a pantalon decomposition if and
only if $r$ is maximal (namely, $X_s(M,P)=\emptyset$ for
$s>r$ and $X_r(M,P)\neq\emptyset$). This maximal $r$ is
given by the formula $r=3g+m+q-3$, and the number of
pantalons is $2g+m+q-2$, where $g$ is the genus of $M$,
$m=|P|$, and $q$ is the number of boundary components of
$M$.

The proof of the following proposition is left to the
reader.

\bigskip\noindent
{\bf Proposition 5.7.} {\sl Assume neither $(M,P)$ is as in
Example 1 or 2, nor $(M,P)=(T^2,\emptyset)$. Let $g$
denote the genus of $M$, $m=|P|$, and $q$ the number of
boundary components of $M$.

\smallskip
(i) Let $A=(a_1,\dots,a_r)$ and $B=(b_1,\dots,b_r)$ be two
generic $r$-families of disjoint curves. The class $[A]$ and the
class $[B]$ are in the same $\MM(M,P)$-orbit if and only
if there exist: 

\smallskip
$\bullet$ a one-to-one correspondence between the
components of $M_A$ and the components of $M_B$, 

\smallskip
$\bullet$ a
permutation $\sigma\in\Sigma_r$, 

\smallskip\noindent
such that every pair
$(N,N')$ composed by a component $N$ of $M_A$ and its
corresponding component $N'$ of $M_B$ satisfies the
following properties:

\smallskip
$\bullet$ the genus of $N$ is equal to the genus of $N'$, the
number of boundary components of $N$ is equal to the
number of boundary components of $N'$, and $|P\cap
N|=|P\cap N'|$;

\smallskip
$\bullet$ if $c:S^1\to\partial N$ is an exterior boundary curve of
$N$, then there exists an exterior boundary curve
$c':S^1\to\partial N'$ of $N'$ such that $\rho_A\circ
c=\rho_B\circ c'$;

\smallskip
$\bullet$ if $a_i$ is a separating limit curve of $N$, then
$b_{\sigma(i)}$ is a separating limit curve of $N'$;

\smallskip
$\bullet$ if $a_i$ is a nonseparating limit curve of $N$, then
$b_{\sigma(i)}$ is a nonseparating limit curve of $N'$.

\smallskip
(ii) Let $A=(a_1,\dots,a_r)$ be a generic family of disjoint curves.
There exist generic closed curves $a_{r+1},\dots,a_s$
(where $s=3g+m+q-3$) so that
$(a_1,\dots,a_r,a_{r+1},\dots,a_s)$ is a generic family of disjoint
curves which determines a pantalon decomposition.

\smallskip
(iii) Let $s=3g+m+q-3$. Let $t_1,t_2,t_3$ be three
non-negative integers such that $2t_1+t_2=m$ and
$t_1+t_2+t_3=2g+m+q-2$. Let $(\tilde M,\tilde P)$ denote
the (non-connected) surface with punctures which is the
disjoint union of $t_1$ pantalons of type I, $t_2$
pantalons of type II, and $t_3$ pantalons of type III. 
Choose $2s$ distinct boundary curves $c_1,\dots,
c_s,c_1',\dots,c_s'$ of $\tilde M$ (provided with the
orientation induced by the one of $\tilde M$). 
We have (up to homeomorphism) $M=\tilde M/\sim$, where we
identify the points of the curves $c_i$ and $c_i'$ for all
$i\in\{1,\dots,s\}$ by $c_i(z)=(c_i')^{-1}(z)$. Denote the
natural projection by $\rho:\tilde M\to M$. Then
$P=\rho(\tilde P)$.
Write $a_i=\rho\circ c_i$ for all
$i\in\{1,\dots, s\}$. Then $A=(a_1,\dots,a_s)$ is a
pantalon decomposition of $(M,P)$, and every 
$\MM(M,P)$-orbit of pantalon decompositions has a
representative which can be obtained in this manner.} \qed

\bigskip
Observe that the above proposition shows that there are 
only finitely
many $\MM(M,P)$-orbits of pantalon
decompositions, and that every generic
family of disjoint curves is a face of a pantalon decomposition, so:

\bigskip\noindent
{\bf Corollary 5.8.} {\sl There are only finitely many
$\MM(M,P)$-orbits in $X(M,P)$.} \qed

\bigskip
We turn now to see that the action of $\MM(M,P)$ on
$X(M,P)$ fournishes examples of pairs of self-commensurating
subgroups, one included in the other, and that this
phenomenon shows the limit of a possible generalization of
Theorem 3.2.

Choose a generic 2-family $(a,b)$ of disjoint curves such that $[a]$
and $[b]$ are in different $\MM(M,P)$-orbits; for example,
take $a=a_1$ and $b=a_4$ in Example 4. Then
$$
\Stab([a,b])\subset\Stab([a])\ .
$$
So, we have the following.

\bigskip\noindent
{\bf Proposition 5.9.} {\sl There exists an example of a
group $G$ and two subgroups $H_0,H_1$ of $G$ such that:

\smallskip
(a) $\Com_G(H_0)=H_0$ and $\Com_G(H_1)=H_1$;

\smallskip
(b) $H_0\subset H_1$;

\smallskip
(c) $H_0$ is not conjugate to $H_1$.} \qed

\bigskip
Assume we are in the situation of Proposition 5.9.
Choose a finite dimensional irreducible unitary
representation $\pi_0$ of $H_0$. The induced representation
$\pi_1=\Ind_{H_0}^{H_1}\pi_0$ is an irreducible representation
of $H_1$, since $\Com_G(H_0)=H_0$ implies $\Com_{H_1}(H_0)
=H_0$. Now, the equality $\Ind_{H_1}^G\Ind_{H_0}^{H_1}=
\Ind_{H_0}^G$ implies:
$$
\Ind_{H_0}^G\pi_0=\Ind_{H_1}^G\pi_1\ .
$$
This shows that every representation of $G$ induced
by a finite dimensional irreducible unitary representation of
$H_0$ is also induced by an (infinite dimensional)
irreducible unitary representation of $H_1$.

\bigskip\bigskip
\centerline{\soustitre 6. The structure of the stabilizers}

\bigskip
In order to fully apply Corollary 5.6, we would like to be
able to exhibit finite dimensional irreducible unitary
representations of the stabilizers, different from the
trivial ones. As pointed out in the introduction, there is
no systematic procedure to produce (finite dimensional)
irreducible unitary representations of a given group, exept for
virtually abelian groups. We turn now precisely
to show that $\Stab([A])$ is virtually abelian if
$A=(a_1,\dots, a_r)$ is a pantalon decomposition.

Define the {\it pure mapping class group} $\PMM(M,P)$ as
the subgroup of $\MM(M,P)$ of isotopy classes of
homeomorphisms $h\in\HH(M,P)$ which pointwise fixe $P$.
Let $\Sigma_P$ denote the group of permutations of $P$. We
have the exact sequence
$$
1\to\PMM(M,P)\to\MM(M,P)\to\Sigma_P\to 1\ .
$$
Notice that the pure mapping class group of a pantalon (of
any type) is the free abelian group generated by the Dehn
twists along its boundary components.

Define the {\it cubic group} $\Cub_r$ to be the group of
linear transformations $f\in GL(\BR^r)$ such that
$f(e_i)\in\{\pm e_1,\dots,\pm e_r\}$ for all
$i\in\{1,\dots,r\}$, where $\{e_1,\dots,e_r\}$ denotes the
canonical basis of $\BR^r$. There is a natural
homomorphism $F_A:\Stab([A])\to \Cub_r$ defined as
follows. Let $\xi\in\Stab([A])$, and let $h\in\HH(M,P)$
which represents $\xi$. Then we set
$$
F_A(\xi)(e_i)=\left\{\matrix{ 
e_j\quad\hfill&{\rm if}\ h(a_i)\simeq a_j\ ,\hfill\cr
-e_j\quad\hfill&{\rm if}\ h(a_i)\simeq a_j^{-1}\ .
\hfill\cr}
\right.
$$
This homomorphism is not surjective in general. Its kernel
is composed by the isotopy classes of homeomorphisms
$h\in\HH(M,P)$ such that $h(a_i)\simeq a_i$ for all
$i\in\{1,\dots, r\}$.

Let $c_1,\dots,c_q$ be the boundary components of $M$.
For $\alpha=[a_1,\dots,a_r]\in X(M,P)$, we denote by
$\TT(\alpha)$ the subgroup of $\MM(M,P)$ generated by
$\tau_{a_1},\dots,\tau_{a_r},
\tau_{c_1},\dots,\tau_{c_q}$.
Recall from Proposition 4.8
that this group is a free abelian group of rank $r+q$.
Notice also that $\TT(\alpha)\subset\Stab(\alpha)$.

\bigskip\noindent
{\bf Proposition 6.1.} {\sl Let $A=(a_1,\dots,a_r)$ be a
generic family of disjoint curves which determines a pantalon
decomposition. Then
$$
\TT([A])=\PMM(M,P)\cap\Ker F_A\ .
$$}

\noindent
{\bf Remark.} Proposition 6.1 implies that $\TT([A])$ has
finite index in $\Stab([A])$ since both
$\PMM(M,P)\cap\Stab([A])$ and $\Ker F_A$ have finite index
in $\Stab([A])$.

\bigskip\noindent
{\bf Proof.} 
The inclusion $\TT([A])\subset\PMM(M,P)\cap\Ker F_A$ being
obvious, it suffices to prove the reverse inclusion:
$\PMM(M,P)\cap\Ker F_A\subset\TT([A])$.
Let $\xi\in\PMM(M,P)\cap\Ker F_A$, and let
$h\in\HH(M,P)$ which represents $\xi$. We have $h\circ
a_i\simeq a_i$ for all $i\in\{1,\dots,r\}$. By Proposition
4.3, we may assume that $h\circ a_i=a_i$ for all
$i\in\{1,\dots,r\}$. By the structure of the pure mapping
class groups of the pantalons, it follows that $\xi$ can
be expressed as
$$
\xi=\tau_{a_1}^{n_1}\dots \tau_{a_r}^{n_r}
\tau_{c_1}^{m_1}\dots\tau_{c_q}^{m_q}\ ,
$$
thus $\xi\in\TT([A])$. \qed

\bigskip
We turn now to show that a given
stabilizer can be expressed by means of mapping class
groups of ``smaller'' surfaces.

Let $A=(a_1,\dots,a_r)$ be a generic family of disjoint curves.
Consider the (non-connected) surface $M_A$ of Section 5
which is the natural compactification of
$M\setminus(\cup_{i=1}^r a_i)$, and denote by $N_1,\dots,
N_l$ the components of $M_A$. We write
$$\eqalign{
\MM(M_A,P)&=\MM(N_1,N_1\cap P)\times\dots\times\MM(N_l,
N_l\cap P)\ ,\cr
\PMM(M_A,P)&=\PMM(N_1,N_1\cap P)\times\dots\times\PMM(N_l,
N_l\cap P)\ .\cr}
$$
We have the exact sequence
$$
1\to\PMM(M_A,P)\to\MM(M_A,P)\to\Sigma_{N_1\cap P}\times
\dots\times\Sigma_{N_l\cap P}\to 1\ .
$$
Recall also that $\rho_A:M_A\to M$ denotes the continuous
map induced by the inclusion map $M\setminus(\cup_{i=1}^r
a_i)\to M$. This map induces a homomorphism
$\rho_\ast:\MM(M_A,P)\to\MM(M,P)$. Notice that, by
Proposition 4.3, 
$$
\Im\rho_\ast=\Ker F_A\subset\Stab([A])\ .
$$

\bigskip\noindent
{\bf Lemma 6.2.} {\sl Assume neither $(M,P)$ is as in
Example 1 or 2, nor $(M,P)=(T^2,\emptyset)$. Let
$A=(a_1,\dots,a_r)$ be a generic family of disjoint curves. For
$i\in\{1,\dots,r\}$, let $c_i$ and $c_i'$ denote the
boundary curves of $M_A$ such that $\rho_A\circ
c_i=\rho_A\circ c_i'=a_i$. Then $\Ker \rho_\ast$ is
generated by 
$\{\tau_{c_1}\tau_{c_1'}^{-1},\dots,\tau_{c_r}
\tau_{c_r'}^{-1}\}$ and is a free abelian group of rank
$r$.}

\bigskip\noindent
{\bf Proof.} Let $\TT$ denote the subgroup of $\MM(M_A,P)$
generated by
$\{\tau_{c_1}\tau_{c_1'}^{-1},\dots,\tau_{c_r}
\tau_{c_r'}^{-1}\}$.
The fact that $\TT$ is a free abelian group of rank $r$ is
a direct consequence of Proposition 4.8, and the inclusion
$\TT\subset\Ker\rho_\ast$ is obvious. So, it remains to
prove the inclusion $\Ker\rho_\ast\subset\TT$.

Let $\xi\in\Ker\rho_\ast$. Observe first that we have the
commutative diagram:
$$\matrix{
1&\to&\PMM(M_A,P)&\to&\MM(M_A,P)&\to&\Sigma_{N_1\cap P}\times
\dots\times\Sigma_{N_l\cap P}&\to&1\cr
&&\downarrow&&\downarrow&&\downarrow\cr
1&\to&\PMM(M,P)&\to&\MM(M,P)&\to&\Sigma_P&\to&1\cr}
$$
and that the homomorphism 
$\Sigma_{N_1\cap P}\times\dots\times\Sigma_{N_l\cap P}\to
\Sigma_P$
is injective, thus $\xi\in\PMM(M_A,P)$.

Let $h:M_A\to M_A$ be a homeomorphism which represents
$\xi$. We have $h\circ c_i=c_i$ and $h\circ c_i'=c_i'$ for all
$i\in\{1,\dots,r\}$. Let $d_1,\dots, d_q$ be the boundary
curves of $M$, and let $e_1,\dots,e_q$ be the corresponding
boundary curves of $M_A$ (i.e. $e_i$ satisfies $\rho_A\circ
e_i=d_i$). We also have $h\circ e_i=e_i$. We choose
essential closed curves $a_{r+1},\dots,a_s$ such that
$(a_1,\dots,a_r,a_{r+1},\dots,a_s)$ is a pantalon
decomposition of $(M,P)$. For $j\in\{r+1,\dots,s\}$, let
$c_j$ denote the essential closed curve
of $M_A$ such that $\rho_A\circ c_j=a_j$. Notice that
$(c_{r+1},\dots, c_s)$ determines a pantalon decomposition
of $(M_A,P)$, and
$c_1,\dots,c_r,c_1',\dots,c_r',e_1,\dots,e_q$ are the
boundary components of $M_A$. Let $j\in\{r+1,\dots,s\}$.
Since $\xi\in\Ker\rho_\ast$, $\rho_A\circ h\circ
c_j\simeq\rho_A\circ c_j$. By Proposition 4.4, it follows
that $h\circ c_j\simeq c_j$. So, by Proposition 4.3, we
may assume $h\circ c_j=c_j$. By the structure of the pure
mapping class groups of the pantalons, it follows that
$\xi$ can be written
$$
\xi=\tau_{c_1}^{u_1}\dots\tau_{c_r}^{u_r}\tau_{c_1'}^{v_1}
\dots\tau_{c_r'}^{v_r}\tau_{c_{r+1}}^{u_{r+1}}\dots
\tau_{c_s}^{u_s}\tau_{e_1}^{w_1}\dots\tau_{e_q}^{w_q}\ ,
$$
where $u_1,\dots,u_s,v_1,\dots,v_r,w_1,\dots,w_q\in\BZ$.
Finally, the equality
$$
1=\rho_\ast(\xi)=\tau_{a_1}^{u_1+v_1}\dots\tau_{a_r}^{u_r
+v_r}\tau_{a_{r+1}}^{u_{r+1}}\dots\tau_{a_s}^{u_s}
\tau_{d_1}^{w_1}\dots\tau_{d_q}^{w_q}
$$
implies by Proposition 4.8:
$$
u_1+v_1=\dots=u_r+v_r=u_{r+1}=\dots=u_s=w_1=\dots=w_q=0\ ,
$$
thus
$$
\xi=(\tau_{c_1}\tau_{c_1'}^{-1})^{u_1}\dots(\tau_{c_r}
\tau_{c_r'}^{-1})^{u_r}\ .\quad\square
$$

\bigskip
The above considerations lead to the following proposition.

\bigskip\noindent
{\bf Proposition 6.3.} {\sl Assume neither $(M,P)$ is as in
Example 1 or 2, nor $(M,P)=(T^2,\emptyset)$. Let
$A=(a_1,\dots,a_r)$ be a generic family of disjoint curves. Then we
have the exact sequence
$$
1\to\BZ^r\to\MM(M_A,P)\to\Stab([A])\to\Cub_r\ .\quad\square
$$}

\bigskip\bigskip
\centerline{\soustitre References}

\bigskip
\item{[BLM]}
J. Birman, A. Lubotzky, J. McCarthy,
{\it Abelian and solvable subgroups of the mapping
class group},
Duke Math. J. {\bf 50} (1983), 1107--1120.

\medskip
\item{[BH]}
M. Burger, P. de la Harpe,
{\it Constructing irreducible representations of discrete
groups},
Proc. Indian Acad. Sci., Math. Sci. {\bf 107} (1997), 223--235.

\medskip
\item{[FLP]}
A. Fathi, F. Laudenbach, V. Po\'enaru,
``Travaux de Thurston sur les surfaces, S\'eminaire
Orsay'',
Ast\'erisque {\bf 66--67}, 1979.

\medskip
\item{[Ha1]}
J. Harer,
{\it Stability of the homology of the mapping class groups
of orientable surfaces},
Ann. Math. {\bf 121} (1985), 215--249.

\medskip
\item{[Ha2]}
J. Harer,
{\it The virtual cohomological dimension of the mapping
class group of an orientable surface},
Invent. Math. {\bf 84} (1986), 157--176.

\medskip
\item{[Harv]}
W.J. Harvey,
{\it Boundary structure of the modular group},
Riemann surfaces and related topics: Proc.
1978 Stony Brook Conf.,
Ann. Math. Stud. 97, Princeton, 1981.

\medskip
\item{[Iv1]}
N.V. Ivanov,
{\it Complexes of curves and the Teichm\"uller
modular group},
Russ. Math. Surv. {\bf 42} (1987), 55--107.

\medskip
\item{[Iv2]}
N.V. Ivanov,
``Subgroups of Teichm\"uller modular groups'',
Translations of Mathematical Monographs 115,
Amer. Math. Soc., Providence, 1992.

\medskip
\item{[IM]}
N.V. Ivanov, J.D. McCarthy,
{\it On injective homomorphisms between Teichm\"uller
modular groups I},
Invent. Math. {\bf 135} (1999), 425--486.

\medskip
\item{[Mac]}
G.W. Mackey,
``The theory of unitary group representations'',
Chicago Lectures in Mathematics,
The University of Chicago Press, Chicago-London,
1976.

\medskip
\item{[MM1]}
H.A. Masur, Y.N. Minsky,
{\it Geometry of the complex of curves I: hyperbolicity},
preprint.

\medskip
\item{[MM2]}
H.A. Masur, Y.N. Minsky,
{\it Geometry of the complex of curves II: hierarchical
structure}, preprint.

\medskip
\item{[McC]}
J.D. McCarthy,
{\it Automorphisms of surface mapping class groups. A 
recent theorem of N. Ivanov},
Invent. Math. {\bf 84} (1986), 49--71.

\medskip
\item{[Pa]}
L. Paris,
{\it Parabolic subgroups of Artin groups},
J. Algebra {\bf 196} (1997), 369--399.

\medskip
\item{[PR1]}
L. Paris, D. Rolfsen,
{\it Geometric subgroups of surface braid groups},
Ann. Inst. Fourier {\bf 49} (1999), 417--472.

\medskip
\item{[PR2]}
L. Paris, D. Rolfsen,
{\it Geometric subgroups of mapping class groups},
preprint.

\medskip
\item{[Ro]}
D. Rolfsen,
{\it Braid subgroup normalisers, commensurators,
and induced representations},
Invent. Math. {\bf 130} (1997), 575--587.

\medskip
\item{[War]}
G. Warner,
``Harmonic analysis on semi-simple Lie groups I'',
Springer -Verlag, Berlin-Heidelberg-New York, 1972.

\bigskip\bigskip
\halign{#\hfill\cr
Luis Paris\cr
Laboratoire de Topologie\cr
UMR 5584 du CNRS\cr
Universit\'e de Bourgogne\cr
BP 47870\cr
21078 Dijon Cedex\cr
FRANCE\cr
lparis@u-bourgogne.fr\cr}

\end